\documentclass{amsart}
\usepackage{graphs,epsfig,graphicx,amsmath,amsfonts,verbatim}

\newtheorem*{thm*}{Theorem}

\newtheorem{lem}{Lemma}
\newtheorem{prop}{Proposition}
\theoremstyle{definition}

\newtheorem*{defs*}{Definition}
\newtheorem{rem}{Remark}

\newcommand{\Z}{\mathbb{Z}}

\begin{document}
\title{Ergodicity of the adic transformation on the Euler
graph}%
\author{SARAH BAILEY}
\address{Department of Mathematics, CB 3250, Phillips Hall,\\
University of North Carolina,\\ Chapel Hill, NC 27599 USA\\
sebailey@email.unc.edu}
 \author{MICHAEL KEANE}
 \address{Department of
Mathematics, Wesleyan University, \\ 265 Church Street,\\
Middletown, CT 06459 USA\\ mkeane@wesleyan.edu}

\author{KARL
PETERSEN} \address{Department of Mathematics, CB 3250, Phillips
Hall,\\ University of North Carolina, \\Chapel Hill, NC 27599 USA\\
petersen@math.unc.edu}

\author{ IBRAHIM A. SALAMA}
\address{School of Business, North Carolina Central
University,\\P.O. Box 19407,\\ Durham, NC 27707 USA\\
salama@cs.unc.edu}

\maketitle

\begin{abstract}
The Euler graph has vertices labelled $(n,k)$ for $n=0,1,2,...$ and
$k=0,1,...,n$, with $k+1$ edges from $(n,k)$ to $(n+1,k)$ and
$n-k+1$ edges from $(n,k)$ to $(n+1,k+1)$. The number of paths from
(0,0) to $(n,k)$ is the Eulerian number $A(n,k)$, the number of
permutations of 1,2,...,$n+1$ with exactly $n-k$ falls and $k$
rises. We prove that the adic (Bratteli-Vershik) transformation on
the space of infinite paths in this graph is ergodic with respect to
the symmetric measure.
\end{abstract}

\bibliographystyle{amsplain}

\section{The Euler Graph}

  The Euler graph is an infinite directed graph such that at level
  $n$ there are $n+1$ vertices labelled $(n,0)$ through $(n,n)$.  The
  vertex $(n,k)$ has $n+2$ total edges leaving it, with $k+1$ edges
  connecting it to vertex $(n+1,k)$ and $n-k+1$ edges connecting
  it to vertex $(n+1,k+1)$.

\begin{center}
\begin{figure}[h]
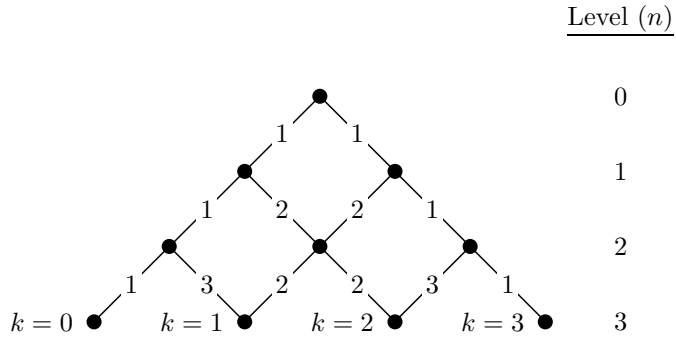

\begin{graph}(4,5)(-2,-2.25)
\roundnode{v00}(-1,1) \freetext(3,2){\underline{Level $(n)$}}
\roundnode{v10}(-2,0)\roundnode{v11}(0,0) \freetext(3,1){0}
\freetext(3,0){1}\roundnode{v20}(-3,-1)\roundnode{v21}(-1,-1)\roundnode{v22}(1,-1)\freetext(3,-1){2}
\edge{v00}{v10}\edge{v00}{v11}\edge{v10}{v20}\edge{v10}{v21}\edge{v11}{v21}\edge{v11}{v22}\edgetext{v00}{v10}{1}
\edgetext{v10}{v20}{1}\edgetext{v00}{v11}{1}\edgetext{v11}{v22}{1}\edgetext{v10}{v21}{2}\edgetext{v11}{v21}{2}
\roundnode{v30}(-4,-2)\roundnode{v31}(-2,-2)\roundnode{v32}(0,-2)\roundnode{v33}(2,-2)\freetext(3,-2){3}\edge{v20}{v30}\
\edge{v20}{v31}\edge{v21}{v31}\edge{v21}{v32}\edge{v22}{v32}\edge{v22}{v33}\edgetext{v20}{v30}{1}\edgetext{v20}{v31}{3}
\edgetext{v21}{v31}{2}\edgetext{v21}{v32}{2}\edgetext{v22}{v32}{3}\edgetext{v22}{v33}{1}\nodetext{v30}(-.7,0){$k=0$}
\nodetext{v31}(-.7,0){$k=1$}\nodetext{v32}(-.7,0){$k=2$}\nodetext{v33}(-.7,0){$k=3$}
\end{graph}
\caption{The first three levels of the Euler graph.  The numbers on
the diagonals give the number of edges coming out of each vertex,
and $k$ represents the label on each vertex.}
\end{figure}
\end{center}

Define $X$ to be the space of infinite edge paths on the Euler
graph.  $X$ is a compact metric space in a natural way: if two paths
$x=x_0x_1x_2...$ and $y=y_0y_1y_2...$ agree for all $n$ less than
$j$ and $x_j\neq y_j,$ then define $d(x,y)=2^{-j}$. The number of
paths from the root vertex (0,0) to the vertex $(n,k)$ is the
Eulerian number, $A(n,k)$, which is the number of permutations of
$0,1,...,n$ with exactly $k$ rises and $n-k$ falls. These numbers
satisfy the recursion
\begin{equation}
A(n+1,k)=(n-k+2)A(n,k-1)+(k+1)A(n,k). \label{recursion}
\end{equation}

\begin{figure}[h]
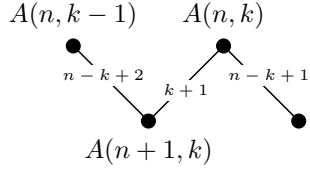

\begin{center}
\begin{graph}(3,2)(0,-.5)
\roundnode{nk1}(0,1)\roundnode{n1k}(1,0)\roundnode{nk}(2,1)\roundnode{n1k1}(3,0)
\edge{nk1}{n1k}\edge{nk}{n1k}\edge{nk}{n1k1}\freetext(.4,.6){\tiny
$n-k+2$}\freetext(1.5,.4){\tiny $k+1$} \freetext(2.6,.6){\tiny
$n-k+1$}\nodetext{nk1}(0,.4){$A(n,k-1)$}\nodetext{nk}(0,.4){$A(n,k)$}
\nodetext{n1k}(0,-.4){$A(n+1,k)$}
\end{graph}
\end{center}

\caption{The Euler Graph gives rise to Equation \ref{recursion}. }
\label{A}
\end{figure}

 We put a partial order on the set of paths in $X$.
The edges $e_0$ through $e_{n+1}$ into the fixed vertex $(n,k)$ with
$0<k<n$ are completely ordered; we illustrate it so that the
ordering increases from left to right. If $x,y$ are paths in $X$, we
say that $x$ is less than $y$ if there exists an $N$ such that both
$x$ and $y$ pass through vertex $(N+1,k)$, $x_n=y_n$ for all $n>N$,
and $x_N <y_N$ with respect to the edge ordering.

\begin{figure}[h]
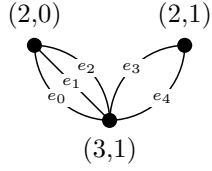

\label{order}
\begin{center}
\begin{graph}(3,2)(0,-.5)
\roundnode{v20}(1,1)\roundnode{v21}(3,1)\roundnode{v31}(2,0)
\nodetext{v20}(0,.4){(2,0)}\nodetext{v21}(0,.4){(2,1)}\nodetext{v31}(0,-.4){(3,1)}
\bow{v20}{v31}{-.2}\edge{v20}{v31}\bow{v20}{v31}{.2}\bowtext{v20}{v31}{-.2}{\tiny
$e_0$}\edgetext{v20}{v31}{\tiny $e_1$}\bowtext{v20}{v31}{.2}{\tiny
$e_2$}\bow{v21}{v31}{-.2}\bow{v21}{v31}{.2}\bowtext{v21}{v31}{-.2}{\tiny
$e_3$}\bowtext{v21}{v31}{.2}{\tiny $e_4$}
\end{graph}
\end{center}
\caption{The order on the edges coming into vertex (3,1):
$e_0<e_1<e_2<e_3<e_4$.}
\end{figure}

Define $k_n:X\to \{0,1,\dots,n-1\}$ by agreeing that if a path $x$
passes through vertex $(n,k)$, then $k_n(x)=k$.  We then say that
$x$ has a \emph{left turn} at level $n$ if $k_{n+1}(x)=k_n(x)$ and a
\emph{right turn} if $k_{n+1}(x)=k_n(x)+1$. Then define
$X_{\text{max}}$ to be the set of paths in $X$ such that there are
no greater paths with respect to the above ordering:
 $X_{\text{max}}$= \{ the path with no left turns, the path with no right
turns\}$\cup$\{$x\in X|$there is a $j$ such that $x$ has a unique
left turn at $x_j$ and for all $n\geq j$, $x_n$ is the maximal edge
into $(n,j+1)$\}. $X_{\text{min}}$ is the set of paths in $X$ such
that there are no smaller paths with respect to the above ordering:
 $X_{\text{min}}$= \{ the path with no left turns, the path with no right
turns\}$\cup$\{$x\in X|$ there is a $j$ such that $x$ has a unique
right turn at $x_j$ and for all $n\geq j$, $x_n$ is the minimal edge
into $(n,n-j)$\} Both $X_{\text{max}}$ and $X_{\text{min}}$ are
countable.

\begin{center}
\begin{figure}[h]
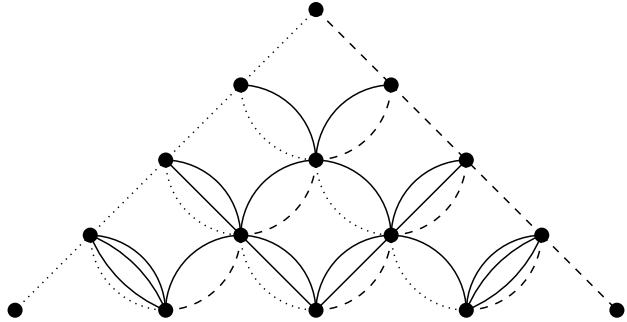

\label{maxmin}
\begin{graph}(4,5.5)(1,-1.25)
\roundnode{-11}(-1,-1)\roundnode{1-1}(1,-1)\roundnode{3-1}(3,-1)\roundnode{5-1}(5,-1)
\roundnode{7-1}(7,-1)
\roundnode{00}(0,0)\roundnode{20}(2,0)\roundnode{40}(4,0)\roundnode{60}(6,0)
\roundnode{11}(1,1)\roundnode{31}(3,1)\roundnode{51}(5,1)\roundnode{22}(2,2)
\roundnode{42}(4,2)\roundnode{33}(3,3)
\edge{22}{00}[\graphlinedash{.5 2}]\edge{42}{60}[\graphlinedash{ 3
3}]\bow{22}{31}{-.2}[\graphlinedash{.5 2}]
\bow{42}{31}{.2}[\graphlinedash{ 3
3}]\bow{11}{20}{.2}\edge{11}{20}\bow{11}{20}{-.2}[\graphlinedash{.5
2}]\bow{31}{20}{-.2}
\bow{31}{40}{.2}\bow{31}{40}{-.2}[\graphlinedash{.5
2}]\bow{51}{40}{-.2}\edge{51}{40}\bow{51}{40}{.2}[\graphlinedash{ 3
3}] \edge{33}{22}[\graphlinedash{.5
2}]\bow{22}{31}{.2}\bow{31}{20}{.2}[\graphlinedash{ 3 3}]
\edge{33}{42}[\graphlinedash{ 3 3}]\bow{42}{31}{-.2}
\edge{00}{-11}[\graphlinedash{.5
2}]\bow{00}{1-1}{-.2}[\graphlinedash{.5
2}]\bow{00}{1-1}{.2}\bow{00}{1-1}{.1}\bow{00}{1-1}{-.1}
\bow{20}{1-1}{.2}[\graphlinedash{ 3
3}]\bow{20}{1-1}{-.2}\edge{20}{3-1}\bow{20}{3-1}{.2}\bow{20}{3-1}{-.2}[\graphlinedash{.5
2}] \edge{40}{3-1}\bow{40}{3-1}{.2}[\graphlinedash{ 3
3}]\bow{40}{3-1}{-.2}
\bow{40}{5-1}{.2}\bow{40}{5-1}{-.2}[\graphlinedash{.5 2}]
\bow{60}{5-1}{-.2}\bow{60}{5-1}{.2}[\graphlinedash{ 3
3}]\bow{60}{5-1}{.1}\bow{60}{5-1}{-.1}
\edge{60}{7-1}[\graphlinedash{ 3 3}]
\end{graph}
\caption{The dashed paths are maximal, and the dotted paths are
minimal.  In addition, the paths following the far left edge and the
far right edge are both maximal and minimal.}
\end{figure}
\end{center}

  If $x\in X\setminus
X_{\text{max}},$ consider the first non-maximal edge, $x_j$, of $x$
and let $y_j$ be the next greatest edge with respect to the edge
ordering.  Then define $y_0y_1...y_{j-1}$ to be the minimal path
into the source of $y_j$ and let
$T(x)=y_0...y_jx_{j+1}x_{j+2}...$(so $T(x)_i=x_i$ for all
$i=j+1,j+2,\dots$). Then $T:X\setminus X_{\text{max}}\to X\setminus
X_{\text{min}}$ is the \emph{Euler adic}.

\begin{center}
\begin{figure}[h]
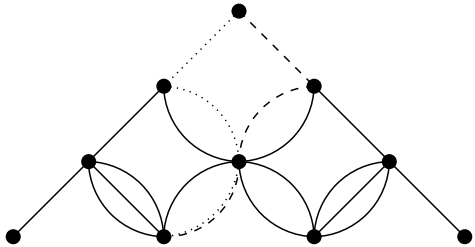

\label{map}
\begin{graph}(4,3.5)(1,-.25)
\roundnode{00}(0,0)\roundnode{20}(2,0)\roundnode{40}(4,0)\roundnode{60}(6,0)
\roundnode{11}(1,1)\roundnode{31}(3,1)\roundnode{51}(5,1)\roundnode{22}(2,2)
\roundnode{42}(4,2)\roundnode{33}(3,3)\edge{22}{00}\edge{42}{60}\bow{22}{31}{-.2}
\bow{42}{31}{.2}\bow{11}{20}{.2}\edge{11}{20}\bow{11}{20}{-.2}\bow{31}{20}{-.2}
\bow{31}{40}{.2}\bow{31}{40}{-.2}\bow{51}{40}{-.2}\edge{51}{40}\bow{51}{40}{.2}
\edge{33}{22}[\graphlinedash{.5
2}]\bow{22}{31}{.2}[\graphlinedash{.5 2}]\bow{31}{20}{.185}
[\graphlinedash{.5 2}]\edge{33}{42}[\graphlinedash{3
3}]\bow{42}{31}{-.2}[\graphlinedash{3 3}]
\bow{31}{20}{.215}[\graphlinedash{3 3}]
\end{graph}
\caption{$T$ maps the dotted path into the dashed path.}
\end{figure}
\end{center}

  Since both $X_{\text{max}}$ and $X_{\text{min}}$ are
  countable, for any $T$-invariant, nonatomic measure $\mu$,
  $\mu(X_{\text{max}})=\mu(X_{\text{min}})=0$.

\section{The Symmetric Invariant Measure}
A \emph{cylinder set} $C=[c_0c_1...c_{n-1}]$ is $\{x\in
X|x_i=c_i\text{ for all } i=0,1,...,n-1\}$.  Given any $T-$invariant
Borel measure, $\mu$, on $X$, define the \emph{weight} $w_n$ on an
edge $c_n$ connecting level $n$ and $n+1$ to be
$\mu([c_0...c_n]|[c_0...c_{n-1}])$ for $n$ greater than 0 and
$w_0=\mu([c_0])$. Then $\mu[c_0...c_n]=w_0...w_n$, where $w_i$ is
the weight on the edge $c_i$.  There are two conditions which
together are necessary and sufficient to ensure that a measure on
$X$ is $T$-invariant. The first is that if $e_0$ and $e_1$ have the
same source vertex and the same terminal vertex, then their weights
are equal.  The second is the \emph{diamond law}. If $u_1$ is the
weight associated with the edges connecting vertex $(n,k)$ to
$(n+1,k)$, $u_2$ is the weight associated to the edges connecting
$(n,k)$ to $(n+1,k+1)$, $v_1$ is the weight on the edges connecting
$(n+1,k)$ to $(n+2,k+1)$, and $v_2$ is the weight on the edges
connecting $(n+1,k+1)$ to $(n+2,k+1)$, then $u_1v_1=u_2v_2$.

\begin{figure}[h]
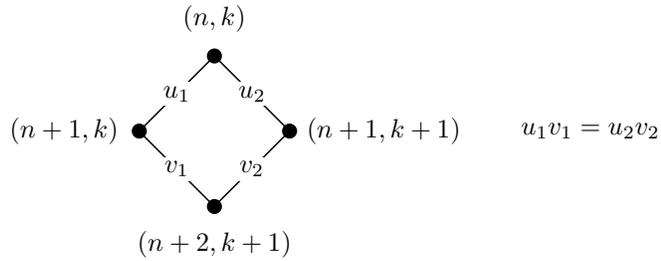

  \begin{graph}(3,3.5)(0,-.75)
  \roundnode{10}(1,0)\roundnode{01}(0,1)\roundnode{12}(1,2)\roundnode{21}(2,1)
  \edge{10}{01}\edge{12}{01}\edge{12}{21}\edge{21}{10}\edgetext{10}{01}{$v_1$}
  \edgetext{12}{01}{$u_1$}\edgetext{12}{21}{$u_2$}\edgetext{21}{10}{$v_2$}
  \freetext(6,1){$u_1v_1=u_2v_2$}\nodetext{10}(0,-.5){$(n+2,k+1)$}\nodetext{01}(-1,0){$(n+1,k)$}
\nodetext{21}(1.25,0){$(n+1,k+1)$}\nodetext{12}(0,.5){$(n,k)$}
  \end{graph}
  \caption{The diamond law.}
  \label{diamond}
  \end{figure}

\begin{defs*}\label{EulerMeasure}
The \emph{symmetric measure}, $\eta$, is determined by assigning
weights 1/$(n+2)$ on each edge connecting level $n$ to level $n+1$.
\end{defs*}

This measure clearly satisfies both of the above conditions and
hence is $T$-invariant.

\begin{figure}[h]
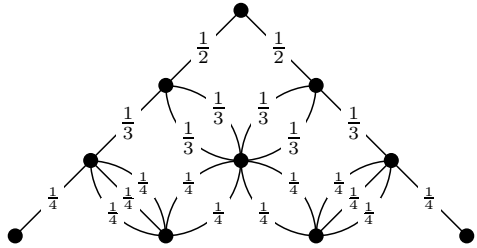

\begin{graph}(6,4)(0,-.25)
\roundnode{00}(0,0)\roundnode{20}(2,0)\roundnode{40}(4,0)\roundnode{60}(6,0)
\roundnode{11}(1,1)\roundnode{31}(3,1)\roundnode{51}(5,1)\roundnode{22}(2,2)
\roundnode{42}(4,2)\roundnode{33}(3,3)\edge{33}{00}\edge{33}{60}\bow{22}{31}{-.2}\bow{22}{31}{.2}
\bow{42}{31}{.2}\bow{42}{31}{-.2}\bow{11}{20}{.2}\edge{11}{20}\bow{11}{20}{-.2}\bow{31}{20}{-.2}\bow{31}{20}{.2}
\bow{31}{40}{.2}\bow{31}{40}{-.2}\bow{51}{40}{-.2}\edge{51}{40}\bow{51}{40}{.2}
\edgetext{33}{22}{$\frac{1}{2}$}\edgetext{33}{42}{$\frac{1}{2}$}\edgetext{22}{11}{$\frac{1}{3}$}
\bowtext{22}{31}{.2}{$\frac{1}{3}$}\bowtext{22}{31}{-.2}{$\frac{1}{3}$}\bowtext{42}{31}{.2}{$\frac{1}{3}$}
\bowtext{42}{31}{-.2}{$\frac{1}{3}$}\edgetext{42}{51}{$\frac{1}{3}$}
\edgetext{11}{00}{\tiny $\frac{1}{4}$}\bowtext{11}{20}{-.2}{\tiny
$\frac{1}{4}$}\edgetext{11}{20}{\tiny $\frac{1}{4}$}
\bowtext{11}{20}{.2}{\tiny $\frac{1}{4}$} \bowtext{31}{20}{.2}{\tiny
$\frac{1}{4}$} \bowtext{31}{20}{-.2}{\tiny $\frac{1}{4}$}
\bowtext{31}{40}{.2}{\tiny $\frac{1}{4}$}\bowtext{31}{40}{-.2}{\tiny
$\frac{1}{4}$}\bowtext{51}{40}{.2}{\tiny$\frac{1}{4}$}\edgetext{51}{40}{\tiny$\frac{1}{4}$}
\bowtext{51}{40}{-.2}{\tiny$\frac{1}{4}$}\edgetext{51}{60}{\tiny$\frac{1}{4}$}
\end{graph}
\caption{The Symmetric Measure}
\end{figure}

\section{The Cutting and Stacking Representation}
  We can also view the transformation $T$ as a map on the unit
  interval defined by ``cutting and stacking" which preserves Lebesgue measure, $m$.
  Each stage of cutting and stacking corresponds to a level in the
  Euler graph.  At each stage $n=0,1,2,\dots$ we have $n+1$ stacks
  $S_{n,0},S_{n,1},\dots,S_{n,n}$ (corresponding to the vertices
  $(n,k)$, $0\leq k\leq n$, of the Euler graph).  Stack $S_{n,k}$
  consists of $A(n,k)$ subintervals of $[0,1]$.  Each subinterval
  corresponds to a cylinder set determined by a path of length
  $n$, terminating in vertex $(n,k)$.  The transformation $\tilde
  T$ is defined by mapping each level of the stack, except the
  topmost one, linearly onto the one above it.  This corresponds
  to mapping each non-maximal path of length $n$ to its successor.
   To proceed to the next stage in the cutting and stacking
   construction, each stack $S_{n,k}$ is cut into $n+2$ equal
   substacks.  These are recombined into new stacks in the order
   prescribed by the way $T$ maps their corresponding cylinder
   sets.  In this way, we obtain a Lebesgue measure-preserving
   transformation defined almost everywhere on $[0,1]$.

\begin{figure}[h]

{\centerline{\includegraphics[width= 6cm]{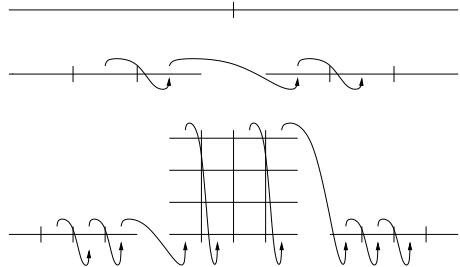}}}

\caption{The Euler adic as a cutting and stacking transformation.}
\end{figure}

\section{Ergodicity}
  In order to prove that the Euler adic $T$ is ergodic with respect to the symmetric measure $\eta$,
we adapt the proof in \cite{Keane} of ergodicity of the
$\mathcal{B}(1/2,1/2)$ measure for the Pascal adic. For previous
proofs of the ergodicity of Bernoulli measures for the Pascal adic,
see \cite{HIK},\cite{vershik1},
\cite{santiago},\cite{mela},\cite{PetMel} and the references that
they contain.

\begin{prop}\label{ldt}For each $x\in X$, denote by $I_n(x)$ the
cylinder set determined by $x_0x_1...x_{n-1}$.  Then for each
measurable $A\subseteq X$,
$$\frac{\eta(A\cap I_n(x))}{\eta(I_n(x))}\to \chi_A(x)\text{  almost everywhere.}$$
\end{prop}

\begin{proof}
In view of the isomorphism of $(X,\eta)$ and $([0,1],m)$, this is
just the Lebesgue Density Theorem.
\end{proof}

Denote by $\rho$ the measure $\eta\times\eta$ on $X\times X$.

\begin{prop}  For $\rho$-almost every $(x,y)\in X\times
X$, there are infinitely many $n$ such that $I_n(x)$ and $I_n(y)$
end in the same vertex of the Euler graph, equivalently
$(n,k_n(x))=(n,k_n(y))$.\label{stacksprop}
\end{prop}

  This is equivalent to saying that for infinitely many $n$ the
  number of left turns in $x_1...x_n$ equals the number of left turns in
  $y_1...y_n$, or that in the cutting and stacking representation
  the subintervals of [0,1] corresponding to $I_n(x)$ and $I_n(y)$
  are in the same stack.  This happens because the symmetric
  measure has a central tendency: if a path is not near the center
  of the graph at level $n$, there is a greater probability that
  at level $n+1$ it will be closer to the center than before (and the farther from the center,
  the greater the probability).  We
  defer momentarily the proof of Proposition \ref{stacksprop} in
  order to show how it immediately implies the main result.

\begin{thm*}  The Euler adic $T$ is ergodic with respect to the
symmetric measure, $\eta$.
\end{thm*}
\begin{proof}
  Suppose that $A\subseteq X$ is measurable and $T$-invariant
  and that $0<\eta(A)<1$.  By Proposition
  \ref{ldt},
  $$\frac{\eta(A\cap I_n(x))}{\eta(I_n(x))}\to 1 \text{ and }
  \frac{\eta(A^c\cap I_n(y))}{\eta(I_n(y))}\to 1\ \text{ for $\rho$-almost
  every } (x,y)\in A\times A^c.$$
Hence for almost every $(x,y)\in A\times A^c$ we can pick an
$n_0=n_0(x,y)$ such that for all $n\geq n_0$,
\begin{equation}
\label{overhalfmeasure} \frac{\eta(A\cap
I_n(x))}{\eta(I_{n}(x))}>\frac{1}{2}\ \ \ \ \ \text{    and    }\ \
\ \ \  \frac{\eta(A^c\cap I_{n}(y))}{\eta(I_{n}(y))}>\frac{1}{2}.
\end{equation}
Then, by Proposition \ref{stacksprop}, we can choose $n\geq n_0$
such that $I_n(x)$ and $I_n(y)$ end in the same vertex, and hence
there is $j\in \Z$ such that $T^j(I_n(x))=I_n(y)$. Since $A$ is
$T$-invariant, this contradicts (\ref{overhalfmeasure}). Then we
must have $\eta(A)=0$ or $\eta(A)=1$, and so $T$ is ergodic with
respect to $\eta$.
\end{proof}

It remains to prove Proposition \ref{stacksprop}.

\begin{lem}  On $(X\times X,\rho)$, for each $n=1,2,\dots$ let
$D_n(x,x')=|k_n(x)-k_n(x')|$, and let
$\mathcal{F}=\mathcal{B}((x_1,x_1'),\dots,(x_n,x_n'))$ denote the
$\sigma$-algebra generated by $(x_1,x_1'),\dots,(x_n,x_n')$.  Let
$\sigma(x,x')$ be a stopping time with respect to $(\mathcal{F}_n)$
such that $D_{\sigma(x,x')}(x,x')>0$.  Fix $M>0$ and let
$$\tau(x,x')=\inf\{n>\sigma(x,x'):D_n\in\{0,M\}\}.$$
For $n=0,1,2,\dots,$ let
\begin{equation*}
Y_n(x,x')=\left\{\begin{array}{lll}
D_{\sigma(x,x')}(x,x')& \text{ if } & 0\leq n\leq \sigma(x,x')\\
D_n(x,x') & \text{ if } & \sigma(x,x')< n \leq \tau(x,x')\\
D_{\tau(x,x')}(x,x') & \text{ if } & n\geq \tau(x,x')\end{array}
\right.
\end{equation*}
Then $(Y_n(x,x'):n=0,1,2,\dots)$ is a supermartingale with respect
to $(\mathcal{F}_n)$.\label{supermartingale}
\end{lem}

\begin{proof}
  We have to check the defining inequality for supermartingales only
  for the range of $n$ where $Y_n=D_n$, since otherwise $Y_n(x,x')$
  is constant in $n$.

  If $x$ turns to the left at stage $n$, then $k_{n+1}(x)=k_n(x)$, but if $x$ turns
  to the right $k_{n+1}(x)=k_n(x)+1$.  From Figure \ref{A} we see
  that
  \begin{equation}\label{prob1}
  \eta\{k_{n+1}(x)=k_n(x)|x_1...x_n\}=\frac{k_n(x)+1}{n+2}
\end{equation}
and \begin{equation}\label{prob2}
\eta\{k_{n+1}(x)=k_n(x)+1|x_1...x_n\}=\frac{n-k_n(x)+1}{n+2}.
\end{equation}
  Without loss of generality assume that $k_n(x')>k_n(x)$.  Note
  that

  \begin{equation*}D_{n+1}= \left\{
     \begin{array}{lll}
      D_n & \text{on the set }& A=\{k_{n+1}(x)=k_n(x),
      k_{n+1}(x')=k_n(x')\}\ \cup\\
      &&\{k_{n+1}(x)=k_n(x)+1,k_{n+1}(x')=k_n(x')+1\} \\
      D_n+1 & \text{on the set } &B=\{k_{n+1}(x)=k_n(x),k_{n+1}(x')=k_n(x)+1\}\\
      D_n-1 &\text{on the set } & C=\{k_{n+1}(x)=k_n(x)+1,k_{n+1}(x')=k_n(x')\}.  \end{array}
     \right.
     \label{Dn}
     \end{equation*}
From (\ref{prob1}) and (\ref{prob2}),

$$
E_{\rho}(D_{n+1}-D_n|\mathcal{F}_n)=0\cdot\rho(A|\mathcal{F}_n)+1\cdot
\rho(B|\mathcal{F}_n)-1\cdot\rho(C|\mathcal{F}_n)$$
$$=\frac{1}{n+2}[(k_n(x)+1)(n-k_n(x')+1)-(k_n(x')+1)(n-k_n(x)-1)]\leq 0.$$

Hence $E_{\rho}(D_{n+1}|\mathcal{F}_n)\leq D_n$.
\end{proof}

\begin{lem}  $\displaystyle \frac{k_n(x)}{n}\to \frac{1}{2}$ in measure.\label{as}
\end{lem}

\begin{proof}
  Let $u_n(x)=2k_n(x)-n$ for all $n$.  We will show that $u_n/n \to 0$
  in measure.  We begin by computing the variance of $u_n$.
  Note that if  $k_{n+1}(x)=k_{n}(x)$ then $u_{n+1}=u_{n}-1,$ and if $k_{n+1}(x)=k_n(x)+1$ then $u_{n+1}=u_{n}+1$.
  Following the calculations in \cite{Salama}, and using (\ref{prob1}) and
  (\ref{prob2}), $$E_\eta(u_{n+1}|u_n)=(n+1)/(n+2)u_n;$$ so, since
  $u_0=0$, $E(u_n)=0$ for all $n=1,2,\dots$ .  Similarly,

  $$\displaystyle E_{\eta}(u_{n+1}^2|u_n)=(u_n-1)^2\left( \frac{k_n(x)+1}{n+2}\right) +
  (u_n+1)^2\left(\frac{n-k_n(x)+1}{n+2}\right)$$

  $$=(u_n-1)^2\left(\frac{u_n+n+2}{2(n+2)}\right)+(u_n+1)^2\left(\frac{n-u_n+2}{2(n+2)}\right)$$
  $$=\frac{nu_n^2}{n+2}+1.$$
  Then   $$E_{\eta}(u_{n+1}^2|u_{n-1})=\frac{n}{n+2}\left(\frac{(n-1)u_{n-1}^2}{n+1}+1\right)+1,$$
  and continuing this recursively we see that
  $$V(u_{n+1})=E_{\eta}(u^2_{n+1})=\frac{1}{(n+1)(n+2)}\sum_{i=0}^{n}(i+1)(i+2)$$
  $$=\frac{n+3}{3}.$$
  Then by Chebyshev's Inequality,
  $$\eta\left\{\left|\frac{u_n}{n}\right|\geq \epsilon\right\}\leq
  \frac{c}{n\epsilon^2}\to 0\text{ as }n\to \infty,$$
  so that $\displaystyle \frac{u_n}{n}\to 0$ in measure, i.e. $\displaystyle \frac{k_n(x)}{n}\to \frac{1}{2}$ in measure.
\end{proof}

\begin{proof}[\textbf{Proof of Proposition 2.}]
  From Lemma \ref{supermartingale}, $(D_n)$ is a supermartingale with respect to
  $\mathcal{F}_n=(\mathcal{B}((x_1,x_1'),\dots,(x_n,x_n')))$.
  Fix $M>0$ and define stopping times $\sigma(x,x')=\text{inf}\{n|k_n(x)\neq
  k_n(x')\}$ and $\tau(x,x')=\text{inf}\{n>
  \sigma(x,x')|D_n\in\{0,M\}\}$.  Then
  $E_{\rho}(D_{\tau}) \leq E_{\rho}(D_\sigma)=1$.  If $\tau$ is
  finite almost everywhere, then
  $$E_{\rho}(D_{\tau})=M(\rho\{D_\tau=M\})+0(\rho \{D_\tau=0\}), \text{ so that }$$
    $\rho\{D_n\neq 0\text{ for any
    }n>\sigma(x,x')\}\leq \rho\{D_\tau=M\}\leq
1/M$ for all $M$.  Letting $M\to\infty$ implies
    that $\rho\{D_n\neq 0\text{ for any }
    n>\sigma(x,x')\}=0$.  Hence with $\rho$-probability 1 there is an
    $n_0$ for which $k_{n_0}(x)=k_{n_0}(x')$.  Repeat this process with $\sigma(x,x')=\text{inf}\{n>n_0(x,x')|
    k_n(x)\neq k_n(x')\}$
    to see that with $\rho$-probability 1, $k_n(x)=k_n(x')$ infinitely many times.  It
    remains to show that $\tau$ is finite almost everywhere.

    We have a fixed $M$; fix also a large $L$.  Fix a small enough
    $\epsilon$ so that if $k_n(x)/n,\ k_n(x')/n$ are in the
    interval $(1/2-\epsilon,1/2+\epsilon)$, then
    $$\frac{k_{n+i}(x)}{n+i}, \frac{n-k_{n+i}(x)}{n+i},
    \frac{k_{n+i}(x')}{n+i},\frac{n-k_{n+i}(x')}{n+i}\geq \frac{1}{4}\text{ for
    }i=0,1,\dots,ML.$$
    In other words, starting from $(n,k_n(x))$ all the
    probabilities of going left or right for both $x$ and $x'$ are
    at least 1/4 for $ML$ steps.  $\text{Let }A_n=\{(x,x')\in X\times
    X|k_n(x)/n,k_n(x')/n\in
    (1/2-\epsilon,1/2+\epsilon)\},$ and note that $\rho(A_n)\to 1$ as $n\to\infty$, by the
    convergence in measure.

    Let $B_n=\{(x,x')\in X\times
    X|k_{n+i}(x)=k_n(x),k_{n+i}(x')=k_n(x')+i\text{ for all }
    i=0,1\dots,M\}$.

     For every $n$, $\{x|\tau(x)=\infty\}\cap A_n\subset A_n\cap
    B_n^c\cap B_{n+M}^c\cap\dots\cap B_{n+(L-1)M}^c=G_n$, since
    $(x,x')$ in $B_n$ implies $D_{n+i}(x,x')$ is either 0 or $M$ for some $i\leq M$.
    Conditioned on the set $A_n$, the sets
    $B_n,B_{n+M},\dots,B_{n+(L-1)M}$ are not independent, because at each
    step the probabilities of going left or right, given by
    sums of the weights on the edges, are changing.  But since the probabilities of going left or right at each step
    are all near 1/2, so that the probability of each event we are considering is
    near the probability that it would be assigned by a genuine symmetric random
    walk, we can estimate the measure of $G_n$.

    For each $j=0,1,\dots,L-1, $ abbreviate
    $E_j=B_{n+jM}^c$.  Then for each pair of vertices
    $v=((jM-1,k),(jM-1,k'))$,we have
    $\rho(E_j|v)\leq (1-1/4^{2M}).$ Thus
    $$\rho(E_j|E_{j-1}\cap\dots\cap E_0\cap
    A_n)=\hskip -.3in\sum_{\begin{array}{l}
    \small{\text{vertices $v$ at}}\\
    \small{\text{level $jM-1$}}\end{array}}\hskip -.3in\rho(E_j|v)\rho(v|E_{j-1}\cap\dots\cap E_0\cap
    A_n)\leq (1-1/4^{2M}),$$
    and iterating gives $\rho(E_{L-1}\cap\dots\cap E_0|
    A_n)\leq \left(1-1/4^{2M}\right)^L.$

Therefore $\rho\time\eta(\tau=\infty|A_n)\leq (1-1/4^{2M})^L$ for
all $L$. Letting $n\to\infty$ and then $L\to \infty$, we conclude
that $\rho\{\tau=\infty\}=0.$
    \end{proof}

\begin{rem}
  \upshape In fact $k_n(x)/n\to 1/2$ almost everywhere.  We can
  see this as follows.  Continue to let $u_n(x)=2k_n(x)-n$ as in Lemma
  \ref{as}.  Since $E((n+2)u_{n+1}|(n+1)u_n)=(n+1)u_n,\ S_n=(n+1)u_n$ forms
  a mean-0 martingale.  If $X_n=S_n-S_{n-1}$, then
  the $X_n$ are a martingale difference sequence in $L^2$, thus
  mean 0 and orthogonal.  The variance of $X_n$ is
  $$E(X_n^2)=E(S_n^2)-E(S_{n-1}^2)=\frac{3n^2+5n}{3}.$$
If we let $b_n=n^2$, then $\sum E(X_n^2)/b_n^2<\infty$, so by the
extension to martingales of Kolmogrov's Criterion for the Strong Law
of Large Numbers (see \cite[p.238]{Feller}) $S_n/b_n\to 0$ almost
everywhere, that is to say, $u_n/n\to 0$ almost everywhere.
\end{rem}

\begin{rem}
\upshape It would be interesting to determine further dynamical
properties of this system, such as weak mixing, rigidity,
singularity of the spectrum, and whether the rank is infinite. So
far we can show that the symmetric measure $(\eta)$ is the only
fully supported invariant ergodic measure \cite{BP}, and that
  $(X,T,\eta)$ is totally ergodic and loosely Bernoulli, \cite{B}.
\end{rem}

  We thank the referee for helpful comments.

\bibliography{bibliography}

\end{document}